\documentclass[a4paper,12pt]{article}
\usepackage{graphicx}
 \usepackage[top=2.5cm,bottom=2.5cm,left=2.5cm,right=2.5cm]{geometry}
 \usepackage{cite, amsmath, amssymb, booktabs, lastpage, graphicx}
\usepackage[hidelinks]{hyperref}
\usepackage{setspace}

\usepackage{amsthm}

\newtheorem{prop}{Proposition}
\newtheorem{lemma}{Lemma}

 \title{\vspace*{2.5cm}{\bf The hitting time index: bounds, a sharp threshold, the proof of a conjecture and the disproof of another conjecture} }
\date{}
\author{{\bf Jos\'e Luis Palacios}\\
{\it Electrical and Computer Engineering Department,}\\{\it The University of New Mexico, Albuquerque, NM 87131, USA}\\
{\bf Shunya Tamura}\\
{\it Okegawa City Okegawa West Junior High School}\\{\it 3680-1 Kawataya, Okegawa, 363-0027, Saitama, Japan}\\
{\bf Jason Wiberg} \\
{\it Electrical and Computer Engineering Department,}\\{\it The University of New Mexico, Albuquerque, NM 87131, USA}\\
and\\
{\bf Aleksandar Petojevi\'c}\\{\it Faculty of Education, University of Novi Sad}\\{\it Podgori\v{c}ka 4, 25000 Sombor, Republic of Serbia}}

\begin{document}
 \thispagestyle{empty} 
\maketitle

\begin{abstract}
 We prove some new upper and lower bounds for the hitting time index ($HT$ index) of an $n$-vertex graph, showing along the way that the largest order of the index is $n^5$ and the smallest is $n^3$. More specifically, we prove the conjecture that the smallest value of the index is $\frac{n(n-1)^2}{2}$, which is attained by the complete graph.  We also find a sharp threshold for the minimal degree $\delta$ so that the $HT$ index goes from an $n^4$ order to an $n^3$ order if the threshold is surpassed.  We disprove the conjecture that the largest value of the index is attained by the $n$-path, and look at the lollipop graphs, some of which attain the $n^5$ order, and seem to be the natural candidates among which the largest value of the index should be attained.  
 
 \end{abstract}

\section{Introduction}

A graph $G=(V, E)$ with vertex set $V=\{1, 2, \ldots, n\}$ and edge set $E\subset V\times V$ will be understood throughout this article to be a simple, connected,  undirected graph. The degree $d_i$ of the vertex $i$ is the number of neighbors of $i$, and we will denote by $\delta$ the smallest degree of $G$. For any details regarding graphs, the reader is directed to reference \cite{BH}.

The random walk on $G$ is the Markov chain $\{ X_n, n\ge 0\}$ with state space $V$ and with transition probabilities $P(i, j)$ from $i$ to $j$,  given by $P(i, j )=\frac{1}{d_i}$ in case $(i,j)\in E$ and 0 otherwise. The hitting time of the state $j$ is the random variable defined as
$$T_j=\inf\{n\ge 0 : X_n=j\},$$
and its expectation when the process is started at $i$ is denoted by $E_iT_j$.  The return time to the state $j$ is defined as
$$T^+_j=\inf\{n\ge 1: X_n=j\},$$
and the mean return time of $j$ is denoted by $E_jT^+_j$. We remark that $E_jT_j=0$, whereas $E_jT^+_j\ge 2$.  In fact, $E_jT^+_j=\frac{2|E|}{d_j}$. One can find in reference \cite{DS} a complete discussion of expected hitting times and mean return times of Markov chains.

Graphs are used in Mathematical Chemistry to model molecules, where the atoms are represented by the vertices and the atomic bonds by the edges.  On these structures, one can define topological descriptors, which are real-valued functions defined on the set of all graphs.  Many topological descriptors, also called indices,  are defined in terms of the degrees of the vertices, others resort to a variety of parameters of the vertices (distances between them, effective resistance between them, eccentricities, etc.) while others look at functions of the eigenvalues of some matrix associated to $G$, like the Laplacian matrix, etc.  All these indices attempt to capture some physico-chemical properties of the molecules, and to classify the molecules according to the values of their descriptor.  We defined in \cite{P} a probabilistic descriptor, the hitting time index $HT$ of the graph $G$, defined by
$$
HT(G)=\sum_{i<j}\max\{E_iT_j, E_jT_i\},
$$
and we showed inequalities relating the $HT$ index to the random walk index, the Kirchhoff index and the Wiener index.  Of those, we will be using these inequalities:
\begin{equation}
\label{elec}
|E|R(G)\le HT(G)\le 2|E| R(G),
\end{equation}
showing the close relationship between the $HT$ index and the Kirchhoff index defined as
$$R(G)=\sum_{i<j} R_{ij},$$
where $R_{ij}$ is the effective resistance between $i$ and $j$ when every edge in $G$ is supposed to be a resistor with unit resistance. Also in \cite{P}, we computed in closed form the $HT$ index for some families of graphs and gave the conjecture that the graphs that maximize and minimize $HT$ are the $n$-path $P_n$ and the complete graph $K_n$, respectively. In \cite{H} they studied closely the relationship between the $HT$ index and the Kirchhoff index, gave some formulas for the $HT$ index of bi-stars and broom graphs (the broom graph $B_{n, d}$ is a $d$-path with $n-d$ pendant vertices attached to one end of the path), and found upper and lower bounds for the $HT$ index in terms of the maximum degree, the hyper-Wiener index, the Wiener index and the Harary index.  In particular, they found the following lower bound for trees:
\begin{equation}
\label{pobre}
HT(T)\ge \frac{H(T)}{2}+\frac{W_2(T)}{2},
\end{equation}
where, for any graph $G$,  $H(G)$ is the Harary index, 
$$H(G)=\sum_{i<j}\frac{1}{d(i,j)}$$
and
$$W_2(G)=\sum_{i<j}d(i,j)^2$$
is the generalized Wiener index with parameter 2.

In \cite{O}, \cite{PaPe},  and their references, the reader may find some more bounds and closed form formulas for the $HT$ index of some families of graphs.

In this article we improve the bound (\ref{pobre}). We also find some new bounds depending only on $n$, showing that the minimal and maximal orders of the $HT$ index are in fact $n^3$ and $n^5$. More specifically, we prove the conjecture that the minimum value of the $HT$ index is uniquely attained by the complete graph. We also find a sharp threshold for the minimum degree $\delta$, showing that the $HT$ index has order $n^3$ when the threshold is surpassed, while there are graphs just below the threshold whose $HT$ index has order $n^4$, and finally, we disprove the conjecture that the graph that attains the maximum value of the $HT$ index is the $n$-path $P_n$, and discuss the lollipop graphs, some of which attain the $n^5$ order and should be taking the place of the $n$-path in a reformulated maximality conjecture. Whenever we say that the order of $HT(G)$ is $n^k$ we mean that $HT(G)=\Theta(n^k)$.

\section{New bounds}

We start by improving the bound (\ref{pobre}) with a simple proof. 
We will  be using the following result found in \cite{HF}
	\begin{lemma} 
		For any tree $T$ and two vertices $a, b \in T$, let $P$ be the unique path from $a$ to $b$,  and let $T_v$ be the connected component (subtree) containing $v$ after removing the edges of $P$.  Then 
		\begin{equation}
			\label{cero}
			E_aT_b = d(a,b)^2+2 \sum_{v \in V(P)}|E(T_v)| d(v, b)
		\end{equation}
		where $|E(T_v)|$ is the number of edges in the subtree $T_v$.
	\end{lemma}
	
	Now we prove the following
	
\begin{prop}
For any tree $T$ we have
\begin{equation}
\label{rico}
HT(T)\ge W_2(T).
\end{equation}
 The equality is attained for $P_2$.
\end{prop}
{\bf Proof.}  From (\ref{cero}) we get immediately that $E_aT_b\ge d(a,b)^2$.  Exchanging the roles of $a$ and $b$, and since $d(b, a)=d(a, b)$ we also get $E_bT_a\ge d(a,b)^2$.  Therefore 
$$\max\{E_aT_b, E_bT_a\}\ge d(a,b)^2,$$
and by the definition of the $HT$ index, we are done. $\bullet$
\vskip .2 in

This bound (\ref{rico}) is always better than (\ref{pobre}) because $\displaystyle \frac{1}{d(i,j)}\le d(i,j)^2$.

\vskip .2 in

Next, we want to give bounds in terms of the number of vertices $n$ of the graph.  One simple strategy is to use existing bounds for the Kirchhoff index, and we can get results like the following

\begin{prop}
For any $G$ we have
$$|E|(n-1)\le HT(G)\le \frac{|E|n(n^2-1)}{3}.$$

The lower bound is attained by the complete graph $K_n$.
\end{prop}

{\bf Proof.} For the lower bound, use $R(G)\ge R(K_n)=n-1$ shown in \cite{L}.  For the upper bound, use $R(G)\le R(P_n)=\frac{n(n^2-1)}{6}$ shown in \cite{P2001}.  $\bullet$
\vskip .2 in
The lower bound, at worst,  is quadratic in $n$, and it is not tight in terms of the order $n^3$, as we shall see below, whereas the upper bound, quintic, is tight in its order $n^5$, though the graph that attains it is not $P_n$, as we shall also see below.  There is an abundant literature of upper and lower bounds for the Kirchhoff index, and we will be using the following
result found in \cite{M}, which provides a much better lower bound than that of the previous proposition:
\begin{lemma} 
For any $G$ with $n\ge 2$ we have
\begin{equation}
\label{crucial}
R(G)\ge \frac{n^2(n-1)-2|E|}{2|E|}.
\end{equation}
The equality is attained by $K_n$.
\end{lemma}

In the next proposition we prove the conjecture proposed in \cite{P}, regarding the minimum value of the $HT$ index, showing along the way that this minimum value is cubic in $n$.

\begin{prop}
For any $G$ with $n\ge 2$ we have
\begin{equation}
\label{mini}
HT(G)\ge HT(K_n)+{n \choose 2} - |E|\ge HT(K_n)=\frac{n(n-1)^2}{2}.
\end{equation}
\end{prop}

{\bf Proof.} Using (\ref{elec}) and (\ref{crucial}) we can write
$$HT(G)\ge |E| R(G)\ge \frac{n^2(n-1)-2|E|}{2} = \frac{n(n-1)^2}{2} +\frac{n(n-1)}{2}-|E|$$
$$=HT(K_n)+{n \choose 2} - |E|\ge HT(K_n),$$
and the last inequality is strict if and only if $G$ is not $K_n$.~~~~$\bullet$
\vskip .2 in

 Besides looking at the Kirchhoff index, as above, in order to find lower bounds for $HT(G)$, we can consider probabilistic arguments, as in the next proposition.






\begin{prop} Let $G$ be any $n$-vertex graph and let $d_1\le d_2\le \cdots \le d_n$ be the degrees of the vertices $1, 2,\ldots, n$ in increasing order.  Then

\begin{equation}
\label{lapropia}
HT(G)\ge |E| \sum_{k=1}^{n-1} {{n-k}\over {d_k}}.
\end{equation}
\end{prop}

{\bf Proof.}  A moment's thought gives us that for any $b\neq a$, $E_aT_a^+\le E_aT_b+E_bT_a$, because the left side of this inequality is the average time to return to $a$ in any possible way while the right hand side requires that this return occurs after a visit to $b$.  Exchanging the roles of $a$ and $b$ we also have $E_bT_b^+\le E_aT_b+E_bT_a$, and therefore
\begin{equation}
\label{max}
\max \{ E_aT^+_a, E_bT^+_b \}\le E_aT_b+E_bT_a.
\end{equation}
Now since 
$$E_aT^+_a=\frac{2|E|}{d_a},$$
we can rewrite (\ref{max}) as
$$\frac{2|E|}{ \min\{d_a, d_b\}}\le E_aT_b+E_bT_a.$$

This implies, in particular, that
$$\max \{E_1T_j, E_jT_1\} \ge \frac{1}{2} (E_1T_j+E_jT_1) \ge \frac {|E|}{d_1},$$
for $2\le j \le n$.
 
Likewise,
$$\max \{E_2T_j, E_jT_2\} \ge \frac {|E|}{d_2},$$
for $3\le j \le n$. Etc.

Now, using these inequalities and the definition of $HT(G)$ we get

$$HT(G)=\sum_{k=1}^{n-1}\sum_{j=k+1}^n \max \{E_kT_j, E_jT_k\} \ge  \sum_{k=1}^{n-1}(n-k) {{|E|}\over {d_k}}.~~~\bullet$$ 

We notice that (\ref{lapropia}), in the case of a $d$-regular graph,  becomes
$$HT(G)\ge \frac{nd}{2}\sum_{k=1}^{n-1}{{n-k}\over {d}}=\frac{1}{4}n^2(n-1),$$
thus obtaining a lower bound which is worse than (\ref{mini}) roughly by a factor of 2.   In fact we can prove  that for all graphs with $|E|$ edges, the worst the lower bound (\ref{lapropia}) can be is obtained when the graph is closest to $d$-regular (for example, a tree can never be regular, but the tree closest to a regular graph is the $n$-path $P_n$, which produces the smallest value of (\ref{lapropia}) among trees).

For that purpose, let us consider the function $\Phi:S\subset \mathbb{R}^n\rightarrow \mathbb{R}$ defined over the set $S$ of $n$-tuples $(d_1, d_2, \ldots, d_n)$ such that $1\le d_1\le d_2 \le \ldots \le d_n\le n-1$ and $\sum_{i=1}^n d_i=2|E|$, given by
	$$\Phi(d_1, d_2, \ldots, d_n)=\sum_{k=1}^n\frac{n-k}{d_k}.$$
Then we have the following

\begin{prop}
The function $\Phi$ attains its minimum in $S$ at the point 
$$d_1=d_2 = \cdots = d_n=\frac{2|E|}{n}.$$
\end{prop}

{\bf Proof.} Since $S$ is closed and bounded, and therefore compact, and $\Phi$ is a continuous function on $S$, $\Phi$ attains a minimum in $S$.

We claim that if there is an $n$-tuple ${\bf d}=(d_1, d_2, \ldots, d_n)\in S$ with a pair $d_i, d_{i+1},$ where the inequality is strict, $d_i<d_{i+1}$, then we can find another $n$-tuple  ${\bf d}^{\epsilon}\in S$ such that  $\Phi({\bf d}^{\epsilon}) <\Phi({\bf d})$ holds.  Indeed, since $d_i<d_{i+1}$, we can find $\epsilon >0$ such that $d_i+\epsilon\le d_{i+1}-\epsilon$.  Consider the $n$-tuple ${\bf d}^{\epsilon}=(d_1, \ldots, d_i+\epsilon, d_{i+1}-\epsilon, \ldots, d_n)$,  i.e., all its coordinates are equal to those of {\bf d} except for the indices $i$ and $i+1$.  That ${\bf d}^{\epsilon} \in S$ is clear because, on the one hand,  the sum of its coordinates is still equal to $2|E|$, and also
$$d_{i-1}\le d_i+\epsilon ~({\rm because~} d_{i-1}\le d_i),$$
$$d_{i}+\epsilon\le d_{i+1}-\epsilon ~({\rm because~ of ~the ~~choice~ of~}\epsilon)$$
$$d_{i+1}-\epsilon \le d_{i+2} ~ ({\rm because~} d_{i+1}\le d_{i+2}).$$

The boundary cases $i=1$ or $i=n-1$ are treated similarly, with just two inequalities to be verified.  Now consider the real function 
$$f(\epsilon)=\Phi({\bf d}^{\epsilon})-\Phi({\bf d})=\frac{n-i}{d_i+\epsilon}+\frac{n-i-1}{d_{i+1}-\epsilon}-\frac{n-i}{d_i}-\frac{n-i-1}{d_{i+1}}.$$
  It is obvious that $f(0)=0$ and also
  $$f^{'}(0)=-\frac{n-i}{d_i^2}+\frac{n-i-1}{d_{i+1}^2}.$$ 
  Since $0<d_i<d_{i+1}$ we have $\displaystyle \frac{1}{d_i^2}>\frac{1}{d_{i+1}^2}$.  Also, since $n-i>n-i-1\ge 0$ we have
  $$\frac{n-i}{d_i^2}>\frac{n-i-1}{d_i^2}>\frac{n-i-1}{d_{i+1}^2},$$
  therefore $f^{'}(0)<0$ (the case $i=n-1$ needs a separate treatment with the same conclusion), implying that $f(\epsilon)<0$ for all sufficiently small $\epsilon >0$, and this means that we can make $\Phi({\bf d}^{\epsilon}) <\Phi({\bf d})$ for a certain $\epsilon >0$.  Thus no point in $S$ with a strict inequality $d_i<d_{i+1}$ can be the minimizer.  The only point in $S$ without strict inequalities between any consecutive coordinates is found when
  \begin{equation}
  \label{equalities}d_1=d_2 = \cdots = d_n=\frac{2|E|}{n},
  \end{equation}
  and this unique candidate left,  must be the minimizer. $~~~~~\bullet$
\vskip .2 in




We remark that the proof of proposition 5 considers continuous variables and the conditions in (\ref{equalities}) may not be met for our discrete $d_i$'s, unless the graph is regular.  It is not difficult to see that in its natural discrete setting, the minimizer is the most balanced integer sequence
$$d_1 = d_2 = \cdots d_{n-r} =q, ~~~~~~~~d_{n-r+1} = \cdots = d_n = q+1,$$
where 
$$q=\left\lfloor \frac{2|E|}{n}\right\rfloor, ~~~r=2|E| -qn,$$
but we omit the details.  Even though (\ref{lapropia}) is worse than (\ref{mini}) for regular graphs, it is easy to find examples where (\ref{lapropia}) is better than (\ref{mini}).  For instance, attach two pendant vertices to two different vertices of a $K_4$ in order to get a graph with 6 vertices and 8 edges.  Simple calculations yield that (\ref{mini}) becomes 82 for this graph, whereas (\ref{lapropia}) becomes 87.3333. This means that (\ref{mini}) and (\ref{lapropia}) are non-comparable.
\vskip .2 in
Looking now for upper bounds of the $HT$ index,  it can be shown (see \cite{P2010})  that for $G$ regular we have
$$\max_{i, j} E_iT_j\le 3n^2,$$
and thus, we immediately obtain the following result:

\begin{prop}
For all regular graphs we have
$$HT(G)\le {n \choose 2} 3 n^2.$$
\end{prop}

The order of this upper bound is quartic,   and it is attained by the $n$-cycle, as seen in \cite{P}. 

Let us denote by $diam(G)$ the largest value of the distances $d(i, j)$ for any pair of vertices $i, j\in V$.  Then we get the following

\begin{prop}
For all $G$ we have
\begin{equation}
\label{otra}
HT(G)\le 2|E|R(G)\le 2|E| {n \choose 2} \max_{i,j}\{R_{ij}\}\le 2|E|{n \choose 2} diam(G).
\end{equation}
\end{prop}

{\bf Proof.} The last inequality is the only one worth an argument, and it holds because the effective resistance $R_{ij}$ between the vertices $i$ and $j$ is bounded above by the distance $d(i, j)$ between $i$ and $j$. Indeed, we choose a shortest path between $i$ and $j$ and create a new graph by deleting from the original graph as many edges as necessary so that this shortest path is the only path between $i$ and $j$. In the new graph, the effective resistance is precisely $d(i,j)$, and by Rayleigh's monotonicity law (see \cite{DS}) this is larger than the effective resistance in the original graph.~~~$\bullet$

 In the case of trees, this yields
$$HT(T)\le n(n-1)^2 diam(T),$$
which aligns with the fact (see \cite{O}) that the broom graphs $B_{n, d}$ with $d$ small (small diameter) have cubic $HT$ values (this upper bound plus the lower bound above), while those of the form $B_{n, n-k}$,  with $k$ small, and thus a large diameter, can be quartic. More generally, it proves that all trees have at most order $n^4$.

\vskip .2 in
So far we  have seen examples of the $HT$ index taking values of orders either $n^3$ or $n^4$.  We report now the following sharp threshold phenomenon between these two powers of $n$ for the $HT$ index in the following proposition.
\vskip .2 in
\begin{prop}
For any $G$ for which $\delta \ge \lfloor \frac{n}{2}\rfloor$, the value of $HT(G)$ is cubic in $n$.   On the other hand, there are graphs with $\delta \le \lfloor \frac{n}{2}\rfloor -1$ such that $HT(G)$ is quartic in $n$.
\end{prop}

{\bf Proof}. It is shown in \cite{P2016} that if $\delta \ge \lfloor \frac{n}{2}\rfloor$ then $R(G)\le 3n-1$.  This together with the bound $HT(G)\le 2|E|R(G)$ implies that $HT(G)$ for these graphs is at most cubic in $n$.  But the inequality (\ref{mini}) says that $HT(G)$ must be at least cubic, thus the first half of the proposition is proved.

Also in \cite{P2016}, one can find an example of a graph $G$ which is  $( \frac{n}{2} -1)$-regular (take $n$ to be even, to simplify the discussion) and $R(G)\ge \frac{n^2}{8}$.  For this example, $|E|$ is quadratic in $n$, and this plus the inequality $HT(G)\ge |E|R(G)$ implies that the order of $HT(G)$ is at least quartic, and an application of proposition 6 completes the proof that the order of $HT(G)$ is indeed quartic.~~~~$\bullet$

\section{A disproof of a conjecture}

In this section we will consider the $L(n, k)$ lollipop graph, which consists of a $k$-clique (a complete graph on $k$ vertices) with an $n-k$-path attached through an extra edge to any vertex of the clique.  

It was conjectured in \cite{P} that the maximum value of the $HT$ index among $n$-vertex graphs is attained by the $n$-path $P_n$, and it was mentioned in that reference that the conjecture is true for $n\le 5$.  However, the conjecture is false in general, and the most economical example is the $L(6,4)$ lollipop, for which $HT(L(6,4))=224$, whereas $HT(P_6)=223$.  

In Table 1 we present the lollipop graphs that attain the maximum value of the $HT$ index for several values of $n$.  These results are loosely aligned with the fact (see \cite{BW}) that the largest expected hitting time of a random walk on an $n$-vertex graph, with value roughly $\frac{4n^3}{27}$,  is attained by the $L(n, \frac{2n}{3})$ lollipop.  The largest value of the $HT$ index is not attained by this graph, but rather by the $L(n, .558n)$ lollipop. The results in Table 1 were obtained by comparing the $HT$ index of every lollipop graph  $L(n,k)$ for each value of $n$. We also verified up to $n=10$ that the lollipop graph maximizes the HT index over the class of all non-isomorphic connected simple graphs by generating every such graph with the \texttt{geng} program from the Nauty and Traces package \cite{geng}. For $n>10$ we found the maximum among lollipop graphs only; whether these values are the global maxima remains an open question.  

\begin{table}[ht]
\centering
\caption{Lollipop graphs $L(n,k)$ that attain the greatest $HT$ index.}
\vskip .2 in
\label{tab:lollipop_ht}
\begin{tabular}{ccc}
\toprule
$n$ & Clique size $k$ & $HT$ index \\
\midrule
6   & 4  & 224.0 \\
7   & 4  & 466.5 \\
8   & 5  & 870.8 \\
9   & 5  & 1,526.4 \\
10  & 6  & 2,531.7 \\
15  & 9  & 17,922.7 \\
20  & 11 & 73,541.4 \\
30  & 17 & 544,547.5 \\
50  & 28 & 6,868,253.4 \\
100 & 56 & 216,963,976.4 \\
500 & 279 & 671,813,209,127.3 \\
1000 & 558 & 21,475,047,362,569.9 \\
\bottomrule
\end{tabular}
\end{table}

So the old conjecture must be changed, replacing the $n$-path with what seems to be, asymptotically, the $L(n, .558n)$ lollipop graph.  What the $n$-path seems to attain is the maximum value of the $HT$ index among trees, and there is  evidence supporting that conjecture in the recent articles \cite{O} and \cite{PaPe}, where the $HT$ index of some families of trees is computed, and all the values found are bounded above by $HT(P_n)$.

Lastly, we are going to show that the $HT$ indices of some lollipop graphs attain the $n^5$ order mentioned after proposition 2. We will be using Chandra's formula (see \cite{Ch}) , which shows that
the commute time between vertices $a$ and $b$,  $E_aT_b+E_bT_a$, can be expressed in electric terms as 
\begin{equation}
\label{commute}
E_aT_b+E_bT_a=2|E|R_{ab}.
\end{equation}

Let us consider now the $L(n,  \frac{n}{2})$ lollipop, and let us assume that $n$ is a multiple of 4, so that all necessary parameters are integers. Let $b$  be the contact point between the $K_{\frac{n}{2}}$ clique and the $P_{\frac{n}{2}}$ path, pick $a\neq b$ in the clique, and pick $c=c_{\frac{n}{4}}$ as the midpoint in the $P_{\frac{n}{2}}$ path. (We are labelling the vertices in the path as $c_1, c_2, \ldots, c_{\frac{n}{2}}=d$). We want to compute $E_aT_c$.  By the nature of the graph under consideration, it is clear that 
$$E_aT_c=E_aT_b+E_bT_c.$$
Now, 
$$E_dT_b=\left(\frac{n}{2}\right)^2=E_dT_c+E_cT_b =\left(\frac{n}{4}\right)^2+E_cT_b,$$
and solving for $E_cT_b$ we get
$$E_cT_b=\frac{3}{16}n^2.$$
By Chandra's formula,  we have
$$E_bT_c+E_cT_b=2|E|R_{bc} = 2\left[{\frac{n}{2} \choose 2}+\frac{n}{2}\right] \frac{n}{4}=\frac{n^3}{16}+\frac{n^2}{8},$$
and solving for $E_bT_c$ we get
$$E_bT_c= \frac{n^3}{16}-\frac{n^2}{16}.$$
And then
$$E_aT_c=E_aT_b+E_bT_c = \frac{n^3}{16}-\frac{n^2}{16}+ \frac{n}{2}-1.$$

Now we argue that instead of $a$ we can choose  $\frac{n}{2}-1$ starting vertices, and instead of $c$ we can  choose $\frac{n}{4}$ target vertices (starting from $a$, to hit $c_j$ with $j > \frac{n}{4}$ takes longer than to hit $c$) and therefore, the hitting time index is larger than a cubic term multiplied by the number of choices for $a$ and by the number of choices for $c$, that is,
$$HT\left(L\left(n, \frac{n}{2}\right)\right)\ge \left(\frac{n^3}{16}-\frac{n^2}{16}+\frac{n}{2}-1\right)\left(\frac{n}{2}-1\right)\left(\frac{n}{4}\right).$$

So we see that the $HT$ index of this lollipop graph is bounded below by a quintic term.

\end{document}